\documentclass[12pt]{amsart}
\usepackage{amscd,amssymb}
\theoremstyle{latex 2e}
\newtheorem{thm}[subsection]{Theorem}
\newtheorem{lem}[subsection]{Lemma}
\newtheorem{prop}[subsection]{Proposition}
\newtheorem{cor}[subsection]{Corollary}

\numberwithin{equation}{section}

\newcommand{\Z}{{\mathbb Z}}

\begin{document}

\title[Finite isometry groups of $4$-manifolds with positive sectional curvature] %
{Finite isometry groups of $4$-manifolds with positive sectional curvature}

\author{Fuquan Fang}
\thanks{Supported partially by
NSF Grant 19925104 of China, 973 project of Foundation Science of
China and the Max-Planck Institut f\"ur Mathematik}
\address{Nankai Institute of Mathematics,
Weijin Road 94, Tianjin 300071, P.R.China}
\address{Department of Mathematics, Capital Normal University,
Beijing, P.R.China} \email{ffang@nankai.edu.cn}
\email{fuquanfang@eyou.com}

\begin{abstract}
Let $M$ be an oriented compact positively curved $4$-manifold. Let
$G$ be a finite subgroup of the isometry group of $M$. Among
others, we prove that there is a universal constant $C$ (cf. Corollary
4.3 for the approximate value of $C$), such that if the order of
$G$ is odd and at least $C$, then $G$ is either abelian of rank at
most $2$, or non-abelian and isomorphic to a subgroup of $PU(3)$
with a presentation $\{A, B| A^m=B^n=1, BAB^{-1}=A^r, (n(r-1),
m)=1, r\ne r^3=1(\text{mod }m) \}$. Moreover, $M$ is homeomorphic
to $\Bbb CP^2$ if $G$ is non-abelian, and homeomorphic to $S^4$ or
$\Bbb CP^2$ if $G$ is abelian of rank $2$.
\end{abstract}
\maketitle

\section{Introduction}
\label{sec:intro}

\vskip2mm

\vskip4mm It is one of the most central problem in Riemannian
geometry to study manifolds with positive sectional curvature. In
dimension three the celebrated work of Hamilton \cite{Ha} shows
that spherical space forms are the only manifolds with positive
curvature metrics. The celebrated Hopf problem asks if $S^2\times
S^2$ admits a metric with positive sectional curvature. This
problem remains open. In case the manifold has a continous
symmetry, i.e., it has a non-zero Killing vector field, by
Hsiang-Kleiner \cite{HK}, this manifold is topologically
homeomorphic to $S^4$ or $\Bbb CP^2$ (in the orientable case). An
important quantity of a Riemannian manifold is its isometry group.
Since K. Grove in 1991 proposed to study positively curved
manifolds with large symmetry, considerable advancement has been
accomplished (cf.
\cite{GS}\cite{GZ}\cite{Sh}\cite{GSh}\cite{Wil}\cite{FR}). In
Shankar \cite{Sh}, Grove-Shankar \cite{GSh} and Baza\u{\i}kin
\cite{Ba}, free isometric actions of rank $2$ abelian groups on
positively curved manifolds are found, which answers in negative a
well-known question of S.S.Chern.

In this paper we are addressed to the following question:

{\it Which groups can be realized as the isometry group of a positively curved $4$-manifold?}

Due to the work of \cite{HK} we may restrict our attention to {\it
finite groups} acting on $4$-manifold by isometries. The question
may have a perfect answer once the order of the group is odd and
larger than a certain constant, depending on the Gromov's Betti
number bound \cite{Gr}.

For the sake of simplicity, in the paper we use $\text{Isom}(M)$
to denote the group of orientation preserving isometries of $M$.

\begin{thm}

Let $M$ be an oriented $4$-manifold with positive sectional curvature. Let
$G\subset \text{Isom}(M)$ be a finite subgroup of odd order. Then
there is a universal constant $C$ such that, if the order $|G|\ge
C$, then $G$ is either abelian of rank at most $2$, or non-abelian
and isomorphic to a subgroup of $PU(3)$ with a presentation $\{A,
B| A^m=B^n=1, BAB^{-1}=A^r, (n(r-1), m)=1, r\ne r^3=1(\text{mod
}m) \}$.
\end{thm}

We refer to Corollary 4.3 for a rough estimate of the constant value
$C$. The following result characterizes the topology of the
manifold $M$ in Theorem 1.1 when $G$ is abelian of rank $2$ or
non-abelian.

\begin{thm}

Let $M$ be an oriented $4$-manifold with positive sectional curvature. Let
$G\subset \text{Isom}(M)$ be a finite subgroup of odd order. If
$|G|\ge C$, then

(1.2.1) $M$ is either homeomorphic to $S^4$ or $\Bbb CP^2$,
provided $G$ is abelian of rank $2$.

(1.2.2) $M$ is homeomorphic to $\Bbb CP^2$, provided $G$ is
non-abelian.
\end{thm}

The case of even order is more involved, in this case we can almost characterize the group as a
subgroup of $SO(5)$.

\vskip 3mm

\begin{thm}

Let $M$ be a compact $4$-manifold with positive sectional
curvature. Let $G\subset \text{Isom}(M)$ be a finite subgroup. If
$M$ is neither homeomorphic to $S^4$ nor $\pm \Bbb CP^2$, then
there is a universal constant $C$ such that, if the order $|G|\ge
C$, then $G$ contains an index $2$ subgroup isomorphic to a
subgroup of $SO(5)$.
\end{thm}

We refer to $\S$8 and $\S$9 for a complete list of possible finite
isometry groups of positively curved $4$-manifolds. The following
result estimate the size of the finite isometry groups.

\begin{thm}

Let $M$ be an oriented compact $4$-manifold with positive
sectional curvature. Let $G\subset \text{Isom}(M)$ be a finite
subgroup of order $|G|\ge C$. If $M$ is not homeomorphic to $S^4$,
then $G$ contains a normal cyclic subgroup of index at most $120$.
\end{thm}

The assumption of $M$ is not homeomorphic to $S^4$ is necessary,
since clearly, $\Bbb Z_p\oplus \Bbb Z_p$ acts on $S^4$ by
isometries for any $p$.

By Cheeger the connected sum $\Bbb CP^2\# \Bbb CP^2$ admits a metric with non-negative sectional
curvature. It is interesting to ask how many copies connected sums of $\Bbb CP^2$ can have a metric
with positive sectional curvature. The following result gives an estimate the size of
elementary abelian $2$-group acting on $4$-manifolds with positive sectional curvature and definite
intersection forms.

\vskip 2mm

\begin{thm}

Let $M$ be a compact oriented $4$-manifold with positive sectional curvature and definite
intersection form. If $\Bbb Z_2^k$ acts on $M$ effectively by isometry, then $k\le 4$.
\end{thm}

We believe that one can improve the above estimate to $k\le 3$ if $\chi (M)\ge 3$.
Note that $\Bbb Z_2^3=\langle T_1, T_2, T_3\rangle $ acts on $\Bbb CP^2$ by
isometries with respect
to the Fubini-Study metric, where $T_1, T_2$ acts holomorphically, and $T_3$ acts by complex conjugation
action.

\vskip 2mm

{\bf Acknowledgement:} The work was done during the author's visit
to the Max-Planck Institut f\"ur Mathematik at Bonn. The author
would like to thank the institute for its financial support and
wonderful research atmospheres. The author is very grateful to
Ian Hambleton for invaluable discussions concerning finite group
actions on $4$-manifolds.

\vskip 5mm

\section{Preparations}

\vskip 4mm

A key fact used in the paper is the following classical result due to Frankel.

\begin{thm} [\cite {Fr}]
Let $M$ be a compact manifold of dimension $m$ with positive sectional curvature, and let
$N_1, N_2$ be two totally geodesic submanifolds of dimensions $n_1$ and $n_2$ respectively. If
$n_1+n_2\ge m$, then $N_1\cap N_2$ is not empty.
\end{thm}

Theorem 2.1 implies readily that the fixed point set of any
isometry of a $4$-manifold with positive sectional curvature contains at
most one $2$-dimensional component. For recent development much
generalizing Theorem 2.1, we refer to \cite{Wil}\cite{FMR}.

\begin{thm} [\cite {Gr}]
Let $M$ be a compact manifold of dimension $m$ with non-negative sectional curvature. Then the total Betti
number $\sum b_i(M)\le C(n)$, where $C(n)$ is a constant depending only on $n$.
\end{thm}

The constant $C(n)$ is roughly $10^{10n^4}$.

\vskip 2mm

\begin{thm} [\cite {Mc1}]
Let $M$ be a simply connected compact $4$-manifold. Let $G$ be a finite group acting effectively on $M$
and trivially on homology groups. If the second Betti number $b_2(M)\ge 3$, then $G$ is abelian of
rank at most $2$, and the fixed point set $\text{Fix} (G;M)$ is not empty.
\end{thm}

\vskip 2mm

\begin{thm}[\cite{Wi1}\cite{HLM}]

Any pseudo-free locally linear action by a finite group on a
closed $4$-manifold homeomorphic to a complex projective space is
conjugate to the linear action of a subgroup of $\text{PSU} (3)$
on $\Bbb CP^2$.
\end{thm}

\vskip 2mm

\section{$q$-extent estimates}

The $q$-extent $xt_q(X)$, $q\ge 2$, of a compact metric space $(X, d)$ is, by
definition, given by the following formula:
$$
xt_q(X)={q \choose 2} ^{-1}\text{max} \Bigl\{ \sum _{1\le i<j\le q}
d(x_i, x_j): \{x_i\}_{i=1}^q \subset X \Bigr\}
$$

Given a positive integer $n$ and integers $k, l \in \Bbb Z$ coprime to $n$,
let $L(n; k, l)$ be the $3$-dimensional lens space, the quotient space of a
free isometric $\Bbb Z_n$-action on $S^3$  defined by
$$\psi _{k,l}: \Bbb Z_n \times S^3\to S^3; \text{  } g(z_1, z_2)=
(\omega ^kz_1,\omega ^l z_2 )$$
with $g\in \Bbb Z_n$ a generator, $\omega =e^{i\frac{2\pi}n}$ and
$(z_1, z_2) \in S^3\subset \Bbb C^2$.

Note that $L(n;k, l)$ and $L(n; -k, l)$ (resp. $L(n; l, k)$) are
isometric (cf. \cite{Ya} p.536). Obviously $L(n; -k,l)$ and $L(n;
n-k,l)$ are isometric. Therefore, up to isometry we may {\it
always} assume $k, l\in (0, n/2)$ without loss of generality. The
proof of Lemma 7.4 in \cite{Ya} works identically for $L(n; k, l)$ with
$0<k, l<n/2$ to prove

\vskip2mm

\begin{lem}[\cite{Ya}]
Let $L(n; k, l)$ be a $3$-dimensional lens space of constant sectional
curvature one. Then
$$\aligned
xt_q(L(n;k,l))&\le \text{arccos}\Bigl\{ \text{cos}
(\alpha _q)\text{cos }\pi n^{-\frac 12}\\
&-\frac 12 \{ (\text{cos }\pi n^{-\frac 12}-\text{cos }\pi /n)^2+
\text{sin}^2(\alpha _q) (n^\frac 12 \text{sin }\pi/n -\text{sin }
\pi  n^{-\frac 12})^2 \}^ {\frac 12} \Bigr\}\endaligned
$$
where $\alpha _q=\pi /(2(2-[(q+1)/2]^{-1}))$.
\end{lem}

\vskip2mm

\begin{cor}

Let $L(n;k,l)$ be a $3$-dimensional lens space of constant
sectional curvature one. If $n\ge 61$, then $xt_5(L(n; k, l))<
\pi/3$.
\end{cor}

\vskip2mm

Let $X$ be a compact $4$-manifold with positive sectional
curvature with an effective $G$-action. We call a fixed point of
$G$ is {\it isolated} if $G$ acts freely on the normal unit sphere
of the tangent sphere at the fixed point. The following estimate
is central in the paper.

\vskip2mm

\begin{prop}
Let $X$ be a compact oriented $4$-manifold with positive sectional
curvature. If $\Bbb Z_n$ acts on $X$ by isometries. Then the
action has at most $5$ isolated fixed points.
\end{prop}

\vskip2mm

\begin{proof} We argue by contradiction, assuming $ x_1,...,
x_6$ are six isolated $\Bbb Z_n$-fixed points. Let $\bar
X=X/\Bbb Z_n$. Connecting each pair of
points by a minimal geodesic in $\bar X$, we obtain a
configuration consisting of twenty geodesic triangles. Because
$\bar X$ has positive curvature in the comparison sense,
the sum of the interior angles of each triangle is $>\pi$ and thus
the sum of total angles of the twenty triangles, $\sum
\theta_i>20\pi$. We then estimate the sum of the total angles in
the following way, first estimate from above of the ten angles
around each $\bar x_i$ and then sum up over the six points. We
claim that the sum of angles at $\bar x_i$ is bounded above by
$10\cdot xt_5(\bar x_i)\le 10\frac {\pi}3$ and thus $\sum
\theta_i\le 6(10\cdot \frac{\pi}3)=20\pi$, a contradiction.

Let $S^\perp_{x_i}$ denote the unit
$3$-sphere in the tangent space $T_{x_i} (X)$. Since $\Bbb Z_n$ acts freely
on every $S^\perp _{x_i}$, the open neighborhood of $\bar x_i$ in $\bar X$
is isometric to the metric cone $C(S^3/\Bbb Z_n)$.  By Corollary 3.2
we arrive at a contradiction.
\qed\end{proof}

\vskip 2mm

\section{Lemmas}

\vskip 5mm

\begin{lem}

Let $M$ be a compact $4$-manifold with non-negative sectional curvature. Assume that the second Betti number
$b_2(M)\ge 3$ or $b_2(M)=2$, and $M\ne S^2\times S^2$. If a finite group $G$ acts
effectively on $M$ by isometries, then the normal subgroup $G_0$ of $G$ acting trivially on
$H_2(M)$ trivially satisfies

(i) the index $[G_0:G]\le C$.

(ii) $G_0$ is an abelian group of rank at most $2$.
\end{lem}

\begin{proof} Consider the natural homomorphism $\rho : G\to \text{Aut}(H_2(M))$. By Gromov's theorem,
 $b_2(M)\le C(4)$.
Therefore, the image $\rho (G)$ is a finite subgroup of
$\text{GL}(\Bbb Z, C(4))$. It is well-known that its torsion
subgroup is isomorphic into $GL(\Bbb Z_3, C(4))$ (cf. \cite{Br1})
which clearly has order bounded, a constant depending only on
$C(4)$.

Since $G_0=\ker (\rho)$, by Theorem 2.3 we see that $G_0$ has to be abelian of rank at most $2$.
 The desired result follows.
\end{proof}

\vskip 2mm

It is a standard fact of topology, if $g$ is a self-diffeomorphism
of $M$, the Lefschetz formula reads
 $$\chi (M, g)=\chi (Fix(g))\hspace{2cm} (4.1)
$$
where $\chi (M, g)$ is the Lefschetz number of $g$. In particular, if $g$ acts trivially on homology, i.e.
$\rho (g)=1$, then
$\chi (M)=\chi (Fix(g)) $.

Given a cyclic group $\Bbb Z_m$ acting effectively on $M^4$ with an isolated  fixed point $p\in M$, consider
 the isotropy
representation of $\Bbb Z_m$ at the tangent space $T_pM$. This representation is uniquely determined by a pair of
integers $p, q$, such that the representation is given by $\alpha (u, v)=(e^{\frac{2p\pi i}m}u,e^{\frac{2q\pi i}m}v)$,
 where
$\alpha =e^{\frac{2\pi i}m}\in \Bbb Z_m$ is a generator. Since the
action is effective, the common factor $(p, q)$ must be coprime to
$m$.

\begin{lem}

Let $M$ be a compact $4$-manifold with positive sectional
curvature. If $G$ acts effectively on $M$ by isometries. If
$G_0\subset G$ is the normal subgroup in the above lemma, and it
contains an element of exponent at least $61^4$. Then $\chi(M)\le
7$.
\end{lem}

\begin{proof}

By \cite{Ya} and the above identity (4.1) we may assume that all prime factors of $|G_0|$ is smaller than $61$.
Let $g\in G_0$ be an element of composite order greater than $61$. Consider the cyclic group $\langle g\rangle$.

Case (i) The fixed point set of $\text{Fix}(g)$ contains a surface $\Sigma$.

By Theorem 2.1, it contains at most a $2$-dimensional component.
Consider the fixed points of $g$ outside $\Sigma$, saying, $p_1,
\cdots , p_n$. Note that $\Sigma$ is either a $2$-sphere or a real
projective plane.

If all $p_i$ are isolated fixed points, i.e., $g$ acts freely on
the unit $3$-sphere of the tangent space $T_{p_i}M$ for all $p_i$,
by Proposition 3.3 we know that $n\le 5$. Hence $\chi (M)=\chi
(M^g)\le 7$.

If not, let $p$ be a non-isolated fixed point of $g$ outside
$\Sigma$. Let $g^k$ be an element of $\langle g\rangle$ whose
fixed point set contains a surface $F$ passing through $p$. By
Theorem 2.1 $\Sigma\cap F\ne \emptyset$. Thus, at an intersection
point $q\in \Sigma\cap F$, $g^k$ acts trivially on both the
tangent space $T_q\Sigma $ and $T_qF$. Clearly, $F\ne \Sigma$
implies the tangent spaces $T_qF\ne T_q\Sigma$. Since $g\in G_0$,
$T_qF+T_q\Sigma$ is the tangent space $T_qM$. Thus $g^k$ acts
trivially on $T_qM$. A contradiction since the action is
effective.

Case (ii). The fixed point set $\text{Fix}(g)$ is zero
dimensional.

As in the above, if all fixed points of $g$ are isolated,  by
Proposition 3.3 $\text{Fix}(g)$ contains at most $5$ points.
Therefore by (4.1) we get $\chi (M)\le 5$.

If not, let $p_1, \cdots, p_m$ be the non-isolated fixed points of
$g$. Therefore, there are surfaces $F_1, \cdots, F_m$ passing
through $p_1, \cdots, p_m$ with corresponding isotropy groups of
order $k_1, \cdots , k_m$. We may assume that $k_i$ are all
smaller than $61$, otherwise, we may use the isotropy group
instead of $\langle g\rangle$ to reduce to Case (i). The $F_i$'s
are all orientable surface with positive sectional curvature since $g\in
G_0$, therefore $F_i=S^2$. By Theorem 2.1, $F_i\cap F_j\ne
\emptyset$, for all $1\le i, j\le m$. Observe that the
intersections are zero dimensional if $i\ne j$. Because that $g$
acts on $F_i\cap F_j$, and the action of $g$ on a sphere has only
two isolated points, $F_i\cap F_j$ is contained in
$\text{Fix}(g)$.

Secondly, it is easy to see that no three surfaces, say, $F_1,
F_2, F_3$, can intersect at some point $p_i\in \text{Fix}(g)$.
Therefore, by Theorem 2.1 again, we see that $m\le 3$. Otherwise,
there are two non-intersecting totally geodesic surfaces in $M$.
If $m\le 2$, the desired result follows by Prop. 3.3 and the
formula (4.1).  If $m=3$, first observe that the integers $k_1,
k_2, k_3$ (at most $60$) are pairwisely coprime. Since the order
of $g$ is at least $61^4$, one gets readily a subgroup $G_1\subset
\langle g\rangle$ of order $|G_1|\ge 61$ with only isolated fixed
points. By Proposition 3.3 we get $\chi (M)\le 5$. \end{proof}

By Lemmas 4.1 and 4.2 we get immediately that

\begin{cor}

Let $M$ be a compact $4$-manifold with positive sectional
curvature. Let $G\subset \text{Isom}(M)$. If $|G|\ge C$, then
$\chi(M)\le 7$, where $C=61^8\times |\text{GL}(\Bbb Z_3, C(4))|$,
and $C(4)\approx 10^{10\times 4^4}$ the Gromov constant.
\end{cor}

In the following we will always use $C$ to indicate this universal
constant.

\begin{lem}

Let $M$ be a compact $4$-manifold with positive sectional curvature, and let $G_0$ be as above.
If $b_2(M)\ge3$ or $b_2(M)=2$ but
$M\ne S^2\times S^2$, then $G_0$ is cyclic.
\end{lem}

\begin{proof} If not, there is a subgroup $\Bbb Z_p^2$ in $G_0$.
By \cite{Mc2} we know that $\Bbb Z_p^2$ has exactly $b_2(M)+2$
fixed points. It is easy to see that, at each fixed points, there
are exactly two submanifolds of dimension $2$ with isotropy groups
isomorphic to $\Bbb Z_p$. If $b_2(M)\ge 2$,  there are at least
two totally  geodesic surfaces in $M$ not intersecting to each
other. A contradiction to Theorem 2.1.
\end{proof}

\vskip 2mm

By Lemmas 4.1, 4.2, 4.3, Theorem 2.4 and \cite {Mc3} for $M=S^2\times S^2$,
it is straightforward to verify that

\vskip 2mm

\begin{lem}

Let $M$ be a compact $4$-manifold with positive sectional
curvature. Let $G\subset \text{Isom}(M)$ be a finite group so that
$|G|\ge C$. If $M\ne S^4$, then $G$ contains a normal cyclic
subgroup of index less than $|SL(\Bbb Z_3, 5)|$.
\end{lem}

The assumption of $M\ne S^4$ is necessary, since clearly, $\Bbb
Z_p\oplus \Bbb Z_p$ acts on $S^4$ by isometries for any $p$. The
estimate above is not sharp. Indeed, if $b_2(M)=4, 5$, the fixed
point of some element must contain a surface $\Sigma$, and we can
show much sharper estimate. If $b_2(M)=2, 3$, since the manifold
has only few homeomorphism types, and the automorphisms preserve
the intersection forms, we may have a better estimate.

Recall that a polyhedral group is a subgroup of $SO(3)$. By
Theorems 2.1, 2.4, Lemma 4.4 and \cite{Mc3} it is easy to see that

\vskip 2mm

\begin{prop}
Let $M$ be a compact $4$-manifold with positive sectional
curvature. If $G$ acts on $M$ effectively by isometries, but
trivially on homology, then there exists a constant $C$, such that
if $|G|\ge C$, then either $G$ is cyclic and $3\le b_2(M)\le 5$,
or $b_2(M)\le 2$ and $G$ is a polyhedral group, or a non-cyclic
subgroup of $PU(3)$.
\end{prop}

\vskip 5mm

\section{Odd order isometries}

Let $M$ be an oriented $4$-manifold with positive sectional curvature. Let $G$ be a group of odd order acting on
$M$ by isometries.

\begin{proof}[Proof of Theorem 1.1]

Let $G_0$ be the normal subgroup of $G$ acting homologically trivial on $H_2(M)$. By the first section
we know that $b_2(M)\le 5$, and if $b_2(M)\ge 3$, then $G_0$ is cyclic.
We will prove  $G$ is always abelian of rank at most $2$, in all cases below except case (3).

Case (1) $b_2(M)\ge 3$.

If the fixed point set of $G_0$ contains a surface, $\Sigma$. By the argument before, since $G$ acts on
$\Sigma$, and the order of $G$ is odd, so $\Sigma$ is orientable, and $G$ acts preserving the orientation of
$\Sigma$. Thus $\Sigma =S^2$, and every element of $G/G_0$ acts on $\Sigma$ with (and only) isolated fixed
points. Therefore, by local isotropy representation at a fixed point in $\Sigma$ it is easy to see
that $G_0$ is in the center of $G$. If $G/G_0$ acts non-effectively on
$\Sigma$, it contains a normal subgroup acting trivially on $\Sigma$, which generates a cyclic subgroup
with $G_0$, since both of them acts on the normal two plane of $\Sigma$ freely, and the action of $G$ is
effective on $M$. For the same reasoning above, this enlarged cyclic subgroup is again in the center of
$G$. We keep to use $G_0$ to denote this group. By \cite{Ku} $G/G_0$ is a subgroup of $SO(3)$. Thus
it must be cyclic since the order is odd. This clearly implies that $G$ is an abelian group of rank at
most $2$.

The same argument applies equally to the case when $G_0$ has
non-trivial subgroup with $2$-dimensional fixed point set.
Therefore, we may assume that $G_0$ acts {\it pseudofreely} on
$M$. From the proof of Corollary 4.3 this case we may upgrade the
estimate to $\chi (M)\le 5$. Therefore $b_2(M)=3$, and $G_0$ is
cyclic by Lemma 4.4. Observe that the fixed point set
$\text{Fix}(G_0)$ consists of exactly five points.
 Since $G_0$ is normal subgroup, $G$ acts on the five fixed points, this defines a homomorphism
$\rho : G\to S_5$, the full permutation group of $5$-letters. On the other hand,
the manifold $M$ is homeomorphic to
$\Bbb CP^2\#\Bbb CP^2\#\Bbb CP^2$ or $\Bbb CP^2\# S^2\times S^2$, up to possibly an orientation reversing
(cf. \cite{Fr}). For these manifolds,
the automorphism groups $\text{Aut} (H_2(M))$ has no order $5$ element, indeed, up to $2$-torsion element,
it has order $(2^3-1)(2^3-2)(2^3-4)$. Therefore, the odd order group $\rho (G)\subset S_5$ is isomorphic
to $\Bbb Z_3$, or trivial. In either case, this shows that the action $G$ on the five points has
at least two fixed points. For any such a fixed point, the isotropy representation gives an embedding of
$G$ into $SO(4)$. Because all odd order subgroup of $SO(4)$ is abelian of rank at most $2$,
the desired result follows.

Case (2) $b_2(M)=2$.

In this case, $M$ is homeomorphic to $S^2\times S^2$, $\Bbb CP^2\#\Bbb CP^2$, or $\Bbb CP^2\# \overline
{\Bbb CP^2}$, up to possibly an orientation reversing (cf. \cite {Fr}). It is easy to see that the
automorphisms induced by $G$-action on $M$ on $H_2(M)$ has no nontrivial odd order element. Therefore,
$G=G_0$.

If $G$ contains a subgroup isomorphic to $\Bbb Z_p\oplus \Bbb Z_p$ for some odd prime $p$, by \cite{Mc2}, if
$M$ is not $S^2\times S^2$, then $G$ has four isolated fixed points. Therefore, by isotropy representation
at these points we find two non-intersecting totally geodesic surfaces in $M$. A contradiction to Theorem 2.1.
By \cite{Mc3} $G$ is cyclic, if $M=S^2\times S^2$ and the action is pseudofree. If
$M=S^2\times S^2$ and the action is not pseudofree, we may have two disjoint totally geodesic surfaces among the
isotropy representations of $\Bbb Z_p\oplus \Bbb Z_p$ at the four fixed points. A contradiction to Theorem 2.1
again.  Therefore, we may assume in the following that, {\bf $G$ has no abelian subgroup of rank $2$}
and $M\ne S^2\times S^2$.

For any $g\in G$, by (4.1) the fixed point set has Euler
characterictic $4$. If its fixed point set contains a surface
$\Sigma$, then the rest must be two isolated points (by Theorem
2.1). Therefore, the normalizer $N(g)$ must fix the two points
since its order is odd. This implies that $N(g)$ is a subgroup of
$SO(4)$, and by the assumption we further conclude that $N(g)$
cyclic. The same argument applies also to the case the fixed
points of $g$ all isolated. Therefore, for any $g\in G$, the
normalizer $N(g)$ is cyclic. In particular, this implies that any
$p$-Sylow subgroup of $G$ is cyclic. By the Burnside theorem (cf.
\cite{Wo} page 163) we know that $G$ must contain a cyclic normal
subgroup of finite index, and so $G$ itself is cyclic.

Case (3) $b_2(M)=1$.

Since $H_2(M)\cong \Bbb Z$, the odd order assumption implies that $G$ acts on $H_2(M)$ trivially. By Theorem 2.4
and \cite{Wi2}
the desired result follows.

Case (4) $b_2(M)=0$.

By \cite{Fr} $M\cong S^4$. If the action of $G$ on $M$ is pseudofree, i.e. the singular set consists of
isolated points, by \cite{Ku} we know that $G$ is a polyhedral group or dihedral group. Consequently,
$G$ must be cyclic since $|G|$ is odd. Otherwise, let $g\in G$ is an element with a $2$-dimensional fixed point
set $\Sigma$.
For any $h\in G$, the totally geodesic surface $h(\Sigma)$ intersects with $\Sigma$ at some point $p$. Thus the
isotropy group of $p$ contains the subgroup generated by $g, h$. This subgroup clearly has $p$ as a fixed point,
so it must be a subgroup of $SO(4)$. In particular, $\langle g, h\rangle$ is abelian of rank at most $2$.
This proves that
$g$
is in the center of $G$. Therefore, $G$ acts on $\Sigma$. Obviously the effective part of $G/\langle g
\rangle$ action on $\Sigma$
is cyclic (note that $\Sigma=S^2$). Thus $G$ is an abelian group of rank at most $2$. This completes the proof.
\end{proof}

\begin{proof}[Proof of Theorem 1.2]

(1.2.2) follows readily from the proof of Theorem 1.1. It remains
to prove (1.2.1).

Let $G=\Bbb Z_{pq}\oplus \Bbb Z_p$, where $p$, $q$ are odd
integers. By Corollary 4.3 the second Betti number $b_2(M)\le 5$.
The fixed point set $F$ of $\Bbb Z_{pq}$ consists of the union of
a connected surface $\Sigma $ with at most $5$ points, where
$\Sigma$ may be empty, $S^2$.

If $M=S^2\times S^2$, it is easy to see that $G=G_0$, where $G_0$
is as in Lemma 4.1. The $G$ action has four fixed points.
Therefore, there is a subgroup whose fixed point set is the union
of two surfaces, a contradiction by Theorem 2.1.

Assume now $b_2(M)\ge 3$ or $b_2(M)= 2$ and $M$ has odd type
intersection form. By Theorem 2.3 $G_0$ is cyclic. We may assume
that $G_0=\Bbb Z_{pq}$, and $\Bbb Z_p$ is a subgroup of
$\text{Aut}(H_2(M))$. Since $p$ is odd, by analyzing the
automorphism group it remains only to consider the following
cases:

(i) $b_2(M)=3$, and $M=3\Bbb CP^2$ (or $3\overline{CP}^2$);

(ii) $b_2(M)=5$, and $M=r\Bbb CP^2\#(5-r)\overline{\Bbb CP}^2$, $0\le r\le 5$;

In all these cases, the automorphism group may have order $3$ or
$5$ subgroups, correspondingly $p=3$ or $5$.

In case (i), note that $p=3$ and an order $3$ subgroup of
$\text{Aut}(H_2(M))$ has trace zero. We may assume that the second
factor, $\Bbb Z_3$ of $G$ has trace zero on $H_2(M)$. By the
formula (4.1) the fixed point set of $\Bbb Z_3$ has Euler
characteristic $2$. If $\Bbb Z_{3q}$ has only isolated fixed
points, by (4.1) once again it has exactly $5$ fixed points.
Obviously, $\Bbb Z_3$ keeps at least two points fixed on the five
points. If $\Bbb Z_{3q}$ has a surface $\Sigma =S^2$ in its fixed
point set, $\Bbb Z_3$ acts on the sphere with two isolated fixed
points. By Theorem 2.1 it is easy to see that $G$ can not have
more than three fixed points. Therefore, we may assume in either
cases that $G$ has exactly two fixed points, say $x_1$ and $x_2$.
By the isotropy representations at $x_1$ and $x_2$ there is an
order $3$ element $g\in G$ with a $2$-dimensional fixed point set
so that $tr(g|_{H_2(M)})=0$. By (4.1) again the fixed point set of
$g$ has Euler characteristic $2$ and so $\text{Fix}(g)=S^2$. On
the other hand, by the Atiyah-Singer $G$-signature formula (cf.
\cite{LM} page 266, Corollary 14.9) the $g$-signature
$\text{sig}(g, M)>0$. Since the intersection form of $M$ is
positive definite (or negative definite), by definition
$tr(g|_{H_2(M)})=\text{sig}(g, M)>0$ (or $<0$). A contradiction.

In case (ii), by the proof of Corollary 4.3 we may assume that
$\Bbb Z_{pq}$ has a surface $\Sigma=S^2$ in its fixed point set.
By the Lefschetz fixed point formula (4.1) $\Bbb Z_{pq}$ has
exactly five isolated fixed points outside $\Sigma$. Since $G$
acts on $\Sigma$ with two fixed points (note that $p$ is odd), and
$G$ also acts on the five points, by Theorem 2.1 it remains only
to consider the case of $p=5$ and the second factor of $G$, $\Bbb
Z_5$, acts freely on the five points. Since the automorphism group
$\text{Aut}(H_2(M))$ does not contain $5$-torsion unless $r=0$ or
$5$ (the author thanks Ian Hambleton for pointing out this fact),
and such an order $5$ automorphism has trace zero, the second
factor $\Bbb Z_5$ has trace zero. The same argument in case (i)
applies equally to arrive at a contradiction.
\end{proof}

\vskip 5mm

\section{Proof of Theorem 1.5}

For an involution on compact $4$-manifold $M$, the Atiyah-Singer $G$-signature theorem tells us
$$
\text{sig}(M, T)=\#(F\cap F) \hspace{2cm} (6.1)
$$
where $F$ is the fixed point set of $T$, and $\#(F\cap F)$ is the self-intersection number
of $F$. In particular, if $T$ has no $2$-dimensional fixed point set, then $\text{Sig}(M, T)=0$.

Let $M$ be a compact $4$-manifold with positive definite intersection form. For an involution
as above, by definition we get $\text{sig}(M, T)=\text{tr} (T|_{H^2(M;\Bbb Q)})$. This together
with the formula (4.1) implies that
$$
\chi(\text{Fix}(T))=2+\text{sig}(M, T) \hspace{2cm} (6.2)
$$

\vskip 2mm

\begin{proof}[Proof of Theorem 1.4]

Since $M$ has definite intersection form, it is easy to see that $\Bbb Z_2^k$ acts on the homology group
$H_2(M;\Bbb Z_2)$ by permutation group $S_n$ with respect to the standard basis of $H_2(M;\Bbb Z_2)$,
where $n=b_2(M)$.

For an involution $T_1\in \Bbb Z_2^k$, we claim that $\text{Fix}(T_1)$ is not empty. If not, by
\cite{Br} VII 7.4 we know that $(T_1)_*\in S_n$ is an alternation without fixed letter, otherwise
there is a surface in the fixed point set of $T_1$. Therefore, the trace $\text{tr}T_1|_{H^2(M; \Bbb Q)}=0$.
By (4.1) $\chi (\text{Fix}(T_1))=2$. A contradiction.

If $T_1$ is an
involution with only isolated fixed points, the number of its fixed points is equal to
$2$, otherwise, by (6.2) $\text{sig}(M, T_1)\ne 0$ and then by (6.1) we get a contradiction.
If $T_1$ has only two isolated fixed
points, since $\Bbb Z_2^k$ acts on the set of the two points, there exists a point whose isotropy group contains
subgroup isomorphic to $\Bbb Z_2^{k-1}$. This clearly implies that $T_1$ and $T_2$ has a common fixed
point, for some $T_2\in \Bbb Z_2^k$. Therefore, $T_2$ has a $2$-dimensional fixed point set, $F$,
by the local isotropy representation. Note that $F$ represents a non-trivial homology class in $H_2(M;\Bbb Z_2)$
and the self-intersection number is non-zero, by the definiteness of the intersection form. Consider
$\Bbb Z_2^k$-action on $F$. Let $G_0$ denote the principal isotropy group. We first claim that $G_0$ is generated
by $T_2$, otherwise, there is an element with three dimensional fixed point set. A contradiction since the action
preservs the orientation.
Now $\Bbb Z_2^k/G_0$ acts on $F$, where $F=S^2$ or $\Bbb RP^2$ since the curvature is positive. For the latter it is
easy to see that $\Bbb Z_2^k/G_0$ has $\Bbb Z_2$ rank at most $2$ (cf. \cite{Ku}), and  for the former,
$\Bbb Z_2^k/G_0$ has $\Bbb Z_2$-rank at most $3$, and so $k\le 4$.
\end{proof}

\vskip 5mm

\section{An algebraic theorem}

\vskip 4mm

In this section we consider central extensions of cyclic groups by
polyhedral groups. We will prove that all of these groups can be
embedded in $SO(5)$.

\begin{thm}

Let $G$ be a finite group which is a central extension of a cyclic
group $\Bbb Z_k$ by a group $H$, where $H$ is a finite subgroup of
$SO(3)$. Then $G$ is isomorphic to a subgroup of $SO(5)$.
\end{thm}

Recall that the finite subgroup of $SO(3)$ (i.e. polyhedral group)
is cyclic, a dihedral group, a tetrahedral group $T$, an
octahedral group $O$, or an icosahedral group $I$. It is
well-known that $T\cong A_4$, the alternating group of $4$
letters, $O\cong S_4$, the full permutation group of $4$-letters,
and $I\cong A_5$,  the alternating group of $5$ letters.

By group extension theory, a central extension
$$
1\to \Bbb Z_m\to G\to \bar G\to 1$$ is uniquely characterized by
the $k$-invariant $k\in H^2(\bar G, \Bbb Z_m)$, and the trivial
extension corresponding to the product $\bar G\times \Bbb Z_m$,
has trivial $k$-invariant.

\vskip 3mm

\begin{lem}

The cohomology groups

$H^2(A_4, \Bbb Z_m)\cong \Bbb Z_m\otimes \Bbb Z_6$;

$H^2(S_4, \Bbb Z_m)\cong \Bbb Z_m\otimes \Bbb Z_2$;

$H^2(A_5, \Bbb Z_m)\cong \Bbb Z_m\otimes \Bbb Z_2$;

$H^2(D_{2k},\Bbb Z_m)=0$, if $m$ is odd, and

$H^2(D_{2k},\Bbb Z_m)\cong \Bbb Z_2\otimes \Bbb Z_m$, if $k$ is
odd, and

$H^2(D_{2k},\Bbb Z_m)\cong (\Bbb Z_2\oplus \Bbb Z_2\oplus \Bbb
Z_2)\otimes \Bbb Z_m$, if $k$ is even
\end{lem}

The lemma can be verified via spectral sequence  (cf. \cite{AM}).

\begin{lem}
Let $G$ be a nontrivial central extension of $\Bbb Z_m$ by $S_4$ or $A_5$. Then $G$ is isomorphic to
$\Bbb Z_{m}\otimes _{\Bbb Z_2}O^*$ or  $\Bbb Z_{m}\otimes _{\Bbb Z_2}I^*$, where $O^*$ (resp. $I^*$) is the
binary octahedral group  (resp. binary icosahedral group).
\end{lem}

\begin{proof}

It is clear that the group $\Bbb Z_{m}\otimes _{\Bbb Z_2}O^*$
(resp. $\Bbb Z_{m}\otimes _{\Bbb Z_2}I^*$) is a nontrivial central
extension of $\Bbb Z_m$ by $O$ (resp. $I$). By Lemma 7.2 this is
the only nontrivial central extension, where $m$ is even.
\end{proof}

By Lemma 7.2 it is straightforward to verify the following three
lemmas.

\begin{lem}
Let $G$ be a nontrivial central extension of $\Bbb Z_m$ by $A_4$.

(7.4.1)  If $m=3^rm_+$, $(m_+, 6)=1$, then $G\cong \Bbb
Z_{m_+}\times (\Bbb Z_2\oplus \Bbb Z_2)\rtimes _\alpha\Bbb
Z_{3^r}$ where $\alpha: \Bbb Z_{3^r}\to \Bbb Z_3\to Aut (\Bbb
Z_2\oplus \Bbb Z_2)$ is the unique nontrivial homomorphism;

(7.4.2) If $m=2^rm_+$, $(m_+, 6)=1$, then $G\cong \Bbb
Z_{m_+}\times \Bbb Z_{2^r}\times _{\Bbb Z_2} (Q_8 \rtimes
_\alpha\Bbb Z_{3})$,  where $(Q_8 \rtimes _\alpha\Bbb Z_{3})$ is
the binary tetrahedral group $T^*$.
\end{lem}

If $m=2^r\times 3^sm_+$ and the $k$-invariant of the extension
contains both $2$ and $3$-torsions,

\begin{lem}
Let $G$ be a nontrivial central extension of $\Bbb Z_{2^rm_+\times
3^s}$ by $A_4$ whose $k$-invariant contains both $2$ and
$3$-torsions, where $(m_+,6)=1$. Then $G\cong \Bbb
Z_{2^rm_+}\times _{\Bbb Z_2} (Q_8 \rtimes _\alpha\Bbb
Z_{3^{s+1}})$.
\end{lem}

\begin{lem}
Let $G$ be a nontrivial central extension of $\Bbb Z_{m}$ by
$D_{2k}$, where $D_{2k}$ is the dihedral group of order $2k$ and $k$ is odd.
Then $G\cong \Bbb Z_{m}\times _{\Bbb Z_2}D_{4k}$.
\end{lem}

\begin{proof}[Proof of Theorem 7.1]

Of course $SO(4)$ contains every rank $2$ abelian subgroup. If $G$
is a central extension by a cyclic subgroup of $SO(3)$, it is a
rank $\le 2$ abelian group, and so isomorphic to a subgroup of
$SO(4)$.

Both groups in Lemma 7.3 can be embedded into $SO(4)$ in a
standard way, since $SO(4)=S^3\times _{\Bbb Z_2}S^3$, and $I^*,
O^*$ are subgroups of $S^3$. If the extensions are trivial, both
groups are contained in $SO(3)\times SO(2)\subset SO(5)$.

For a group of type (7.4.1), it is isomorphic to a subgroup of $SO(3)\times SO(2)\subset SO(5)$.
Indeed, it is generated by the subgroup $A_4\subset SO(3)$ and the cyclic subgroup of order
$3^rm_+$ of $SO(3)\times SO(2)$ given by a block diagonal $5\times 5$-matrix whose minor
$3\times 3$ block is the order $3$ permutation  matrix, and the minor $2\times 2$ block is
the order $3^rm_+$ rotation matrix.

For a group in Lemma 7.5 (and also of type (7.4.2)), note that the group $Q_8 \rtimes _\alpha\Bbb
Z_{3^{s+1}}$ is a subgroup of $U(2)$ acting freely on $S^3$ (cf.
\cite{Wo} page 224). Take the cyclic group $\Bbb Z_{2^rm_+}$ in
the center of $U(2)$, it has a common order $2$ subgroup with
$Q_8$ and so it together with $Q_8 \rtimes _\alpha\Bbb
Z_{3^{s+1}}$ generates a subgroup of $U(2)$ of the type in Lemma 7.5.

It remains to consider central extensions by a dihedral group.

By Lemma 7.6 the group in that lemma is clearly a subgroup of $O(4)\subset SO(5)$.  It suffices
to prove the theorem by considering the extension by a dihedral group of
order divisible by $4$. By Lemma 7.2 we may assume further that $m=2^r$.
Consider the $2$-Sylow group of $G$, saying $G_0$, which is a central extension of $\Bbb Z_m$ by
$D_{2^s}$. By Lemma 7.2, the group $H^2(D_{2^s}; \Bbb Z_{2^r})\cong \Z_2^3$, with generators,
$\beta(x), \beta(y), w$, where $\beta$ is the Bockstein homomorphism associated to the extension
$1\to \Bbb Z_2\to \Bbb Z_{2^{r+1}}\to \Bbb Z_{2^{r}}\to 1 $.

If the extension is nontrivial over $\Bbb Z_{2^{s-1}}\subset
D_{2^s}$, then $G$ contains an element generating a subgroup of
index $2$. By the classification of $2$-groups with an index $2$
cyclic subgroup the desired result follows (cf. Wolf's book \cite
{Wo} page 173). Otherwise,  by \cite{AM} page 130 the proof of
Lemma 2.11 we may determine the group $G$. It is easily seen that
the group $G$ is a subgroup of $SO(5)$. This proves the desired
result.
\end{proof}

\vskip 5mm

\section{Even order isometries, I}

In this section we consider group $G$ of isometries of
$4$-manifold $M$ of even order with the assumption,  {\it $G_0$
has no element with $2$-dimensional fixed point set}, where $G_0$
is the normal subgroup of $G$ acting trivially on the homology
group $H_2(M)$. By Lemma 4.4, $G_0$ is cyclic, unless $b_2(M)\le
1$ or $M=S^2\times S^2$. If the order $|G|\ge C$, by Corollary 4.3
$b_2(M)\le 5$, and this may be improved to $b_2(M)\le 3$ under the
above assumption, by Lemma 4.1 and the proof of Lemma 4.2.

\begin{lem}
Let $M$ be a compact oriented $4$-manifold with positive sectional
curvature. Let $G\subset \text{Isom}(M)$ be a finite group. Then
$G$ is an extension of a cyclic group $G_0$ by a subgroup of
$Aut(H_2(M))$. If additionally the intersection form of $M$ is odd
type and $G_0$ has no element with a fixed point set of dimension
$2$, then $G_0$ has odd order.
\end{lem}
\begin{proof} The former one is by Theorem 2.3.
The latter fact is because any involution must have a
$2$-dimensional fixed point, by \cite{Br} VII Lemma 7.6.
\end{proof}

Let $M$ be a compact simply connected $4$-manifold with $2\le
b_2(M)\le 3$. If the intersection form of $M$ is odd type, by
\cite{Fr}, $M$ is homeomorphic to $3\Bbb CP^2$, $2\Bbb CP^2\#
\overline{\Bbb CP}^2$, $2\Bbb CP^2$, and $\Bbb CP^2\#
\overline{\Bbb CP}^2$, by possibly reversing the orientation.

Let $G\subset \text{Isom}(M)$ be a finite group. Consider the
induced homomorphism $\rho : G\to \text{Aut} (H_2(M;\Bbb Z))$. Let
$h: \text{Aut} (H_2(M;\Bbb Z))\to \text{Aut}(H_2(M;\Bbb Z_2))$ be
the forgetful homomorphism.

\begin{lem}
Let $M$ be a compact oriented $4$-manifold with positive sectional
curvature. Let $G\subset \text{Isom}(M)$ be a finite group with
induced homomorphism $\rho : G\to \text{Aut} (H_2(M;\Bbb Z))$. If
$2\le b_2(M)\le 3$ and the intersection form is odd type,  then
$\rho(G)\cap \text{ker}(h)$ can not contain $\Bbb Z_2^2$.
\end{lem}
\begin{proof} First note that the kernel of the forgetful homomorphism $h$
is $\Bbb Z_2^k$, generated by the reflections on the three factors
of $\Bbb Z^k$, where $k=b_2(M)$.  If
 $\rho (G)\cap \text{ker}(h)$ contains $\Bbb Z_2\oplus \Bbb Z_2$, let
 $F$ denote the fixed point set of $\Bbb Z_2\oplus \Bbb Z_2$. By the
 same spectral sequence calculation in \cite{Mc2} with $\Bbb Z_2$-coefficients
 the homology $H^*(F;\Bbb Z_2)$
has rank $b_2(M)+2$, where $F$ has dimension at most $1$.

If $b_2(M)=2$, and $\rho (G)\cap \text{ker}(h)$ contains $\Bbb
Z_2\oplus \Bbb Z_2$, there is an involution $T$ whose trace on
$H_2(M;\Bbb Z)$ is $-2$. By the trace formula (6.1) again the
Euler characteristic $\chi (F(M;T))=0$. A contradiction, since the
fixed point set has positive sectional curvature, and so positive Euler
characteristic.

If $b_2(M)=3$, there is an involution $T_1\in \Bbb Z_2^2$ whose
trace on $H_2(M;\Bbb Z)$ is $-1$ (it can not be $-2$ by (6.1)
since the fixed point set has positive Euler characteristic). By
the trace formula (6.1) the Euler characteristic $\chi
(F(M;T_1))=1$. This together with Lemma 9.1 implies that
$F(M;T_1)=\Bbb RP^2$ (since the positive sectional curvature). For another
involution $T_2\ne T_1\in \Bbb Z_2^2$, since $F$ is equal to the
fixed point set of $T_2$ on $F(M;T_1)$, which is either a point or
$\Bbb RP^1\cup \{pt\}$. Therefore, the homology $H^*(F;\Bbb Z_2)$
has rank at most $3$. A contradiction.
\end{proof}

\vskip 2mm

\begin{lem}
Let $M$ be a compact oriented $4$-manifold with positive sectional
curvature. Let $G\subset \text{Isom}(M)$. Assume that $2\le
b_2(M)\le 3$, and the intersection form of $M$ is odd type. If
$G_0$ has no element with a fixed point set of dimension $2$, then
$G$ is one of the following type

\centerline{ $\Bbb Z_n$, $\Bbb Z_n\rtimes \Bbb Z_2$, $\Bbb
Z_n\rtimes \Bbb Z_4$, $\Bbb Z_n\rtimes (\Bbb Z_2\oplus\Bbb Z_2)$;}

\noindent where $n$ is odd.
\end{lem}

\begin{proof}
There are four cases to consider:

Case (i) If $M=3\Bbb CP^2$;

By Lemma 8.2, $G/G_0$ is a subgroup of a central extension of
$\Bbb Z_2$ by $S_3$, the permutation group of three letters.
Therefore, the $2$-Sylow group of $G$ is $\Bbb Z_2$, $\Bbb Z_2^2$
or $\Bbb Z_4$. By Theorems 1.1 and 1.2 the odd order subgroup of
$G$ is cyclic. Therefore, $G$ is isomorphic to one of the
following:

\centerline{ $\Bbb Z_n$, $\Bbb Z_n\rtimes \Bbb Z_2$, $\Bbb
Z_n\rtimes \Bbb Z_4$, $\Bbb Z_n\rtimes (\Bbb Z_2\oplus\Bbb Z_2)$}

\noindent where $n$ is odd.

Case (ii) If $M=2\Bbb CP^2$;

By Lemma 8.2 again $\rho (G)=G/G_0$ is a subgroup of an extension
of $\Bbb Z_2$ by $S_2=\Bbb Z_2$. It is easy to see that
$G/G_0\cong \Bbb Z_4$ or $\Bbb Z_2$. Hence $G$ is either cyclic,
$G=\Bbb Z_n\rtimes \Bbb Z_2$ or $\Bbb Z_n\rtimes \Bbb Z_4$, where
$n$ is odd.

Case (iii) If  $M=2\Bbb CP^2\# \overline{\Bbb CP}^2$;

As above, by Lemma 8.2 one may show that $G/G_0$ is $\{1\}$, $\Bbb
Z_2$, $\Bbb Z_2\times \Bbb Z_2$, $\Bbb Z_4$. Thus, $G$ is cyclic,
or one of the extensions $\Bbb Z_n\rtimes \Bbb Z_2$, $\Bbb
Z_n\rtimes \Bbb Z_4$, $\Bbb Z_n\rtimes (\Bbb Z_2\oplus\Bbb Z_2)$,
\newline where $n$ is odd.

Case (iv) If  $M=\Bbb CP^2\# \overline{\Bbb CP}^2$;

By Lemma 8.2 it is easy to see that $\rho (G)\cong \Bbb Z_2$. Thus
$G$ is cyclic or $G\cong \Bbb Z_n\rtimes \Bbb Z_2$.
\end{proof}

 \vskip 2mm

 By \cite{Mc2} one sees that $G_0$ is a polyhedral group if $M=S^2\times S^2$. The following lemma is
 essentially due to McCooey.

\begin{lem}
Let $M$ be a $4$-manifold homeomorphic to $S^2\times S^2$. Let $G$
be a finite group acting effectively on $M$ preserving the
orientation. Then, $G$ is an extension of a polyhedral group by a
subgroup of $\Bbb Z_{2}\oplus \Bbb Z_2$. Moreover, if $G$ is an
extension of a polyhedral group by $\Bbb Z_{2}\oplus \Bbb Z_2$,
$G$ is isomorphic to either of

(8.4.1) $D^*_{2^nm_-}\times \Bbb Z_{m_+}$

(8.4.2) $((\Bbb Z_{m_-}\rtimes \Bbb Z_{2^{n+1}})\times \Bbb
Z_{m_+})\rtimes \Bbb Z_2$ with $n>1$;

\noindent  where $m_+, m_-, n$ are all odd integers, $(m_+,
m_-)=1$.
\end{lem}

A complete list of the group $G$ in Lemma 8.4 may be found in
\cite {Mc2}.

\vskip 5mm

\section{Even order isometries, II}

\vskip 5mm

Let $M$ be a compact oriented $4$-manifold with positive sectional
curvature. Let $G\subset \text{Isom}(M)$ be a finite group. Let
$G_0$ be the subgroup of $G$ acting  trivially on homology. In
this section we consider the case that some element of $G_0$ has
$2$-dimensional fixed point set, saying, $\Sigma$. Observe that
$\Sigma= S^2$ or $\Bbb RP^2$, since it is totally geodesic with
positive sectional curvature. Since $G_0$ is normal, $G$ acts on the surface
$\Sigma$. For simplifying the notions, we may assume that $\Sigma$
is a fixed point component of $G_0$ itself (one should keep in
mind it might be a cyclic subgroup of $G_0$ which does the job).
This implies immediately that the effective part of $G/G_0$ action
on $\Sigma$ is a subgroup of $O(3)$ (cf. \cite{Ku}).

\begin{lem}

Let $M$, $G$ and $G_0$ be as above. Assume that $G_0$ is cyclic
with a fixed point surface $\Sigma$. If every element of $G$ has
at least an isolated fixed point in $\Sigma$. Then $G_0$ is in the
center of $G$.
\end{lem}
\begin{proof}
For any $g\in G$ with a fixed point $p\in \Sigma$, since $g$ acts
on $\Sigma$, the isotropy representation of $g$ at $T_pM$ splits
into $T_p\Sigma \oplus \nu _p(\Sigma)$, the tangent and normal
$2$-planes of $\Sigma$ at $p$. Therefore, the group $\langle g,
G_0\rangle$ has a representation of dimension $4$ which splits
into the sum of two dimensional representations. Since $p$ is
isolated, $g$ must preserve the orientation of $T_p\Sigma$ and so
also the orientation of $\nu _p(\Sigma)$. Similarly, $G_0$
preserves the orientation of $ \nu _p(\Sigma)$. This shows that
$G_0$ commutes with $g$.  The desired result follows.
\end{proof}

Case (i) $\Sigma =\Bbb RP^2$.

\begin{lem}

Let $M$, $G$ and $G_0$ be as above. Assume that $b_2(M)\ge 2$, and
$M\ne S^2\times S^2$. If $\Bbb RP^2 \subset \text{Fix}(G_0)$, then
$G$ is a center extension of a cyclic group by either a cyclic
group or a dihedral group.
\end{lem}

\begin{proof}

By Lemma 4.4 $G_0$ is cyclic. Since any isometry of $\Bbb RP^2$
has an isolated fixed point, by Lemma 9.1 $G_0$ is in the center
of $G$. We may assume that $G/G_0$ acts effectively on $\Bbb
RP^2$, otherwise, we may replace enlarge $G_0$ by the extension by
the principal isotropy group of $G/G_0$, and the same argument in
Lemma 9.1 applies equally. Therefore, by appealing to \cite{Ku}
$G/G_0$ is a dihedral group or a cyclic group. The desired result
follows.
\end{proof}

\vskip 3mm

Case (ii) $\Sigma =S^2$.

Since any orientation preserving isometry of $S^2$ has only
isolated fixed points, $G/G_0$ acts pseudofreely on $S^2$, i.e.
the singular data consists of isolated points. Let $G_1\subset
G/G_0$ be the subgroup of isometries preserving the orientation of
$S^2$. Clearly $G_1$ has index at most $2$ in $G/G_0$. By
\cite{Ku} $G_1$ surjects onto a subgroup of $SO(3)$. The kernel of
the surjection acts trivially on $S^2$. Therefore, the same
argument of Lemma 9.1 implies its kernel together with $G_0$
generates a cyclic group in the center of $G$. This proves that

\begin{lem}

Let $M$, $G$ and $G_0$ be as above. Assume that $b_2(M)\ge 2$, and
$M\ne S^2\times S^2$. If $S^2\subset \text{Fix}(G_0)$, then $G$
contains an index $\le 2$ normal subgroup, $H$, which is a central
extension of a cyclic group by a polyhedral group.
\end{lem}

\vskip 3mm

Results in $\S7$ gives a complete classification of the index at
most $2$ subgroup. However, in order to get more precise
information on how the group $G$ might be in Lemma 9.3, we need to
analyze how an orientation reversing element in $G$ conjugates on
$G_0$. There are two types of involution acting on $S^2$ reversing
the orientation, either a reflection with fixed point set a
circle, or free action on $S^2$ with quotient $\Bbb RP^2$. In the
former case we have

\begin{lem}

Let the assumptions be as in Lemma 9.3. If there is an $\alpha \in
G$ so that its fixed point set on $S^2$ has dimension $1$, then
$G$ satisfies an extension
$$
1\to H\to G\to \Bbb Z_2\to 1
$$
where $H$ is as in Lemma 9.3, the conjugation $\alpha g_0\alpha
^{-1}=g_0^{-1}$, for any $g_0\in G_0$.
\end{lem}

\vskip 2mm

\begin{proof}

Let $S=\text{Fix}(\alpha|_{S^2})$. Note that  $\alpha| _{S^2}$ is
a reflection on $S^2$ along $S$. Since $\alpha$ preserves the
orientation of $M$, therefore $\alpha$ is also a reflection with
respect to a line  in the normal plane $\nu _p(S^2)$ (so $\alpha$
has order $2$). Thus, $\alpha$ and $G_0$ generates a dihedral
group. Therefore, the conjugation $\alpha g_0\alpha
^{-1}=g_0^{-1}$, for any $g_0\in G_0$.
\end{proof}

\vskip 2mm

Next let us consider the case that there is an element $\beta \in
G/G_0$ acting freely on $S^2$ (reversing the orientation of
$S^2$). By \cite {Ku} we know that $G/G_0$ is a subgroup of
$O(3)$. The element $\beta $ must be the center of $O(3)$, acting
on $S^2$ by antipodal map.

\begin{lem}

Let the assumptions be as in Lemma 9.3.  If the Euler
characteristic $\chi (M)>2$, and there is an involution $\beta \in
G$ acting freely on $S^2$, then $|G_0|$ is odd.
\end{lem}

\begin{proof}

By Corollary 4.3  $\chi (M)\le 7$. Therefore, if $M$ is of even
intersection type, it is homeomorphic to $S^2\times S^2$, or
$S^2\times S^2\# S^2\times S^2$, by \cite{Fr}. If $|G_0|$ is even,
by a result of Atiyah-Bott \cite{Ab} (compare \cite{LM}), an
involution in $G_0$ must have all components of the same
dimension, since the Spin structure of $M$ is unique. A
contradiction to Theorem 2.1 since $S^2$ is fixed by $G_0$, this
implies all fixed point components (at least two by the trace
formula (6.1)) of $G_0$ are of dimension $2$.

It remains to consider the case when $M$ is non-Spin. If $\chi
(M)$ is odd, by \cite{Br} VII Lemma 7.6 $\beta$ has
$2$-dimensional fixed point set. By Theorem 2.1 this fixed point
set intersects $S^2$, a contradiction to the freeness of $\beta$
on $S^2$. If $\chi (M)$ is even, then $M$ is homeomorphic to
$r\Bbb CP^2\# s\overline{\Bbb CP}^2$, where $r+s=2, 4$. Since
$\beta ([S^2])=-[S^2]\ne 0$ as a homology class, one may verify
case by case that there is always a class $x \in H^2(M;\Bbb Z_2)$,
so that $\beta ^*(x)=x$ and $x^2\ne 0$. By \cite{Br} VII Lemma 7.4
$\beta$ has $2$-dimensional fixed point set again. A
contradiction.
\end{proof}

\begin{lem}

Let the assumptions be as in Lemma 9.5. Then  $G$ satisfies an
extension
$$
1\to H \to G\to G/H\to 1
$$
where $H$ is an extension of the cyclic group $G_0$ by $\Bbb Z_2$,
which is non-splitting when $|G_0|$ is even, and $G/H$ is either a
cyclic or a dihedral group.
\end{lem}

\begin{proof}

Let $p: G\to G/G_0$ be the projection. Let $H=p^{-1}(\langle \beta
\rangle )$. By Lemma 9.5, $H$ is a non-splitting extension of
$G_0$ by $\Bbb Z_2=\langle \beta \rangle$ if $|G_0|$ is even. By
\cite {Ku} $G/G_0$ is a subgroup of $O(3)$. The element $\beta$ is
given by the antipodal map, which belongs to the center of
$G/G_0$. Therefore, the quotient group, $G/H$, acts on
$S^2/\langle \beta\rangle =\Bbb RP^2$. As we noticed in Case (i),
this implies that $G/H$ is either  a cyclic or a dihedral group.
\end{proof}

By Lemma 9.3 we know that the extension in the above lemma
satisfies certain additional commutativity. We leaves it to the
reader.

\vskip 5mm

\vskip 5mm

\section{Proofs of Theorems 1.3 and 1.4}

\vskip 5mm

\begin{proof}[Proof of Theorem 1.3]

Let $G_0$ be the subgroup of $G$ acting trivially on $H_*(M;\Bbb
Z)$. By Corollary 4.3 we may assume that $b_2(M)\le 5$, and
moreover $b_2(M)\le 3$ if $G_0$ contains no element with
$2$-dimensional fixed point set. By Lemma 4.4 we may assume that
$G_0$ is cyclic, or $M=S^2\times S^2$. In this case, by Corollary
4.3, Theorems 8.3 and 8.4 we know that $G$ contains an index $\le
2$ subgroup which is a subgroup of $SO(4)$, since $\Bbb Z_n\rtimes
\Bbb Z_2$ with $n$ odd is always a subgroup of $SO(4)$. On the
other hand, if the action of $G_0$ contains $2$-dimensional
strata, by Lemmas 9.2, 9.3 and Theorem 7.1 we know that $G$
contains an index $\le 2$ subgroup which is a subgroup of $SO(5)$,
provided $M\ne S^2\times S^2$. Finally, if $M=S^2\times S^2$, and
$G_0$-action is not pseudofree, as in the proof of Lemma 9.5,
$|G_0|$ has to be odd. Therefore, by Theorem 1.2 $G_0$ is cyclic
of odd order. Now the same proof of Lemmas 9.2 and 9.3 applies to
conclude that $G$ contains an index $\le 2$ subgroup embedded in
$SO(5)$. The desired result follows.
\end{proof}

\begin{proof}[Proof of Theorem 1.4]

By Lemmas 8.3, 8.4 and Theorem 2.4 we may assume that $G_0$ has
$2$-dimensional strata. By \cite {Wi2} we may assume that
$b_2(M)\ge 2$. If $M=S^2\times S^2$, as in the proof of Theorem
1.3 $G_0$ is cyclic of odd order. The desired result follows since
$\text{Aut}(H_2(M,\Bbb Z))=\Bbb Z_2^2$. Therefore, it remains to
consider the case where $5\ge b_2(M)\ge 2$, and $M\ne S^2\times
S^2$. By Lemmas 9.2 and 9.3, the index $\le 2$ subgroup $H$ is a
central extension of a cyclic group by a polyhedral group.  Since
a polyhedral group other than cyclic or dihedral has order at most
$60$, it suffices to consider when $H$ is a central extension of a
cyclic group by a cyclic group or a dihedral group. Therefore, $H$
contains an abelian subgroup of index $\le 2$ and rank at most
$2$. Therefore, it suffices to prove that this rank $2$ abelian
subgroup contains a cyclic component of index at most $30$.

Assume an abelian group $A$ of rank $2$ with $G_0\subset A$. By
Theorem 1.2 the odd order subgroup of $A$ is cyclic. Clearly $G_0$
has even order. By the proof of Lemma 9.5 we may assume that the
intersection form of $M$ is odd type. Assume that $M=r\Bbb CP^2\#
s\overline{\Bbb CP^2}$, where $r+s\le 5$. It is easy to see that
the automorphism group $\text{Aut}(H_2(M,\Bbb Z))$ is a normal
extension of $\Bbb Z_2^{r+s}$ by a subgroup of $S_{r+s}$, the
permutation group of $(r+s)$ letters. Therefore, the $2$-Sylow
subgroup of $A/G_0$ has order $\le 16$. The desired result
follows.
\end{proof}
\vskip 5mm

%%%%%%%%%%%%%%%%%%%%%%%%%%%

\end{document}